\documentclass[anz]{cupaus}

\newcommand{\pfrac}[2]{\frac{\partial #1}{\partial #2}}
\newcommand{\ddfrac}[2]{\frac{\mathrm{d} #1}{\mathrm{d} #2}}

\newcommand{\e}{\mathrm{e}}

\DeclareMathOperator{\Ei}{Ei}
\DeclareMathOperator{\erf}{erf}
\DeclareMathOperator{\li}{li}

\begin{document}

\verso{Self-similar smoothing of a discontinuity by degenerate cross-diffusion}
\recto{Michael C. Dallaston
}

\title{Self-similar smoothing of a discontinuity by degenerate cross-diffusion}
\author[1]{Michael.~C. Dallaston}
\address[1]{{School of Mathematical Sciences, Queensland University of Technology,  Brisbane, QLD~4000, Australia}
\email{michael.dallaston@qut.edu.au}}

\pages{1}{23}

\begin{abstract}
We consider the dynamics of an idealised cross-diffusion model of biological invasion in which the diffusion of an invading population is inhibited by a resident population, which is in turn degraded by the former.  We formulate and numerically solve the problem that describes the self-similar smoothing of an initially piecewise-constant invading population.  In the limit the height ahead of the initial discontinuity vanishes, the speed of propagation also vanishes, but only logarithmically slowly.  This result is confirmed by a matched asymptotic analysis.  The  singular nature of this limit indicates that the system does not permit the existence of compactly supported solutions that exhibit a moving front or interface, which has important implications for  the simulation of more complex models that also feature this degenerate diffusive mechanism.
\end{abstract}

\keywords[2020 \textit{Mathematics subject classification}]{primary 35Q92, secondary 34E15, 65M08}

\keywords[\textit{Keywords and phrases}]{cross-diffusion, degenerate diffusion, cell invasion, similarity solutions}

\maketitle

\section{Introduction}

\begin{figure}
\centering
\includegraphics{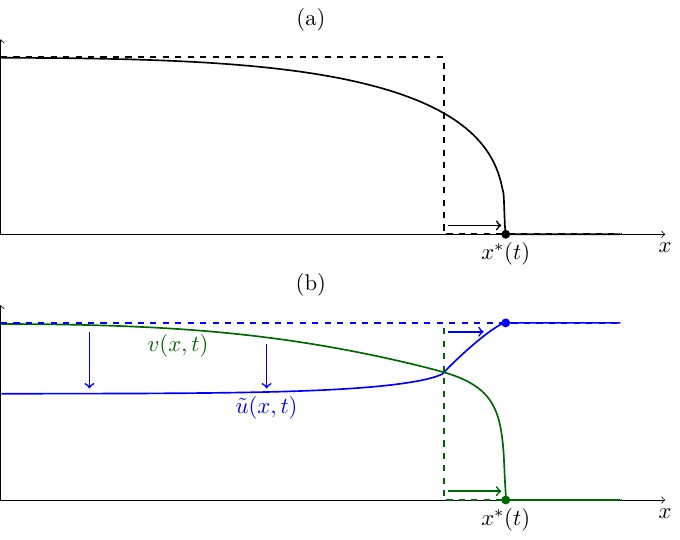}
\caption{(a) Behaviour of the porous medium equation \eqref{eq:pme} from an initially compactly supported solution (dashed line).  The solution over time exhibits a moving interface $x^*(t)$ at which the height $v$ vanishes.  (b) The putative equivalent case for a cross-diffusion model with degradation \eqref{eq:pdeWithReactions}.}
\label{fig:schematic}
\end{figure}

In the field of mathematical modelling of diffusive processes in physical and biological phenomena, some of the most technically challenging and interesting are those that involve degenerate diffusion.  A canonical example of degenerate diffusion is the porous medium equation \cite{aronson2006porous,vazquez2006porous}:
\begin{equation}
\pfrac{v}{t} = \pfrac{}{x}\left(v^m \pfrac{v}{x}\right),
\label{eq:pme}
\end{equation}
where $m$ is a constant exponent.  In this equation, the diffusivity $v^m$ vanishes at points where $v=0$.  As long as $m>0$, given an initial condition that is exactly zero to the right of a point $x^*_0$, for example a piecewise step function
\begin{equation}
v(x,0) = \begin{cases} 1 & x < x^*_0 \\ 0 & x \geq x^*_0 \end{cases},
\label{eq:pmeIC}
\end{equation}
solutions to \eqref{eq:pme} exist that feature a singular point $x^*(t)$; at such a point, the solution (and therefore the diffusivity) vanishes (see Figure~\ref{fig:schematic}(a)).  In the porous media literature such a point is called an interface or free boundary.
Explicit solutions exist that exhibit such an interface, for example, the Barenblatt--Pattle similarity solution \cite{barenblatt1952some,pattle1959diffusion}; an advantage of an initial step function \eqref{eq:pmeIC} is that the resulting solution will also exhibit self-similarity.  Any such solution must be defined in a weak sense, as the derivatives of $v$ in \eqref{eq:pme} will typically not exist at the interface itself.

Adding a logistic reaction term to \eqref{eq:pme} results in the porous Fisher equation, a degenerate-diffusion version of the classical Fisher--Kolmogorov--Petrovsky--Piskunov (Fisher--KPP) problem, and which has seen significant interest in mathematical biology~\cite{aronson1980density,depablo1998travelling,newman1983long,sanchezgarduno1995traveling}.  The degenerate nature of the problem near an interface is unchanged by the inclusion of the reaction term.  In this instance, most interest has been in the existence of travelling wave solutions; the porous Fisher equation has a unique travelling wave solution (in contrast to the continuum of travelling waves of the classical Fisher--KPP problem), and that solution exhibits an interface.

More recently, degenerate reaction-diffusion systems that feature multiple species have been considered in the context of biological cell invasion and interaction.  These models may be complicated by cross-diffusivity, in which the population of one species affects the diffusivity of another.
One such class of two-species models features a non-motile resident cell population $\tilde u$ and a diffusive invading cell population $v$, which are coupled by $\tilde u$ inhibiting the spread of $v$, while at the same time $v$ degrades the local population of $\tilde u$, thus indirectly increasing its own diffusivity.  Nondimensionally, these models may take the form
\begin{subequations}
\begin{align}
\pfrac{\tilde u}{t} &= -\tilde uv + F(\tilde u) \\
\pfrac{v}{t} &= \pfrac{}{x}\left((1-\tilde u)\pfrac{v}{x}\right) + G(\tilde u,v).
\end{align}
\label{eq:pdeWithReactions}%
\end{subequations}
An example of previously studied systems of this form is the approximation \cite{gallay2022propagation,mascia2024numerical} of the Gatenby--Gawlinski model of acid-mediated tumour growth \cite{gatenby1996reaction}, for which $F$ and $G$ are logistic terms in $\tilde u$ and $v$, respectively.  Other similar examples are those considered in \cite{colson2021travelling}, which is similar but neglects the recovery of the resident population $\tilde u$ ($F\equiv0$) and the model considered in \cite{browning2019bayesian,chen2026wavespeed,el2021travelling}, which similarly neglects growth of $\tilde u$ but includes the effects of inhibition of $\tilde u$ on the carrying capacity of $v$, in addition to its diffusivity.  

Models such as \eqref{eq:pdeWithReactions} are simple cases of general Keller-Segel-type  systems that may include both nonlinear, signal-dependent motility (diffusivity) and taxis (that is, flux proportional to gradients of the other species) \cite{fu2012stripe,jin2018boundedness,wang2019boundedness}.   Generally, such systems are qualitatively different to \eqref{eq:pdeWithReactions} (and more regular) as both species diffuse; in such models, $\tilde u$ takes the role of a signalling agent rather than a resident population.  However, systems that include taxis-like `volume-filling' effects in the flux, but a nondiffusive resident population as in \eqref{eq:pdeWithReactions}, have been recently considered in the context of continuum limits of discrete agent-based cell models \cite{crossley2025existence,crossley2023traveling,simpson2024discrete}.

The degeneracy of models of the form \eqref{eq:pdeWithReactions} arises when the resident cell population is at its capacity.  An initial configuration that explores this degeneracy, and is also relevant from the perspective of modelling cell invasion, has the resident population initially at its capacity, with the invading population initially localised to a compactly supported region.  An example, similar to \eqref{eq:pmeIC}, is thus
\begin{equation}
\tilde u(x,0) = 1, \qquad v(x,0) = \begin{cases} 1 & x < x^*_0 \\ 0 & x > x^*_0 \end{cases}.
\label{eq:pdeWithReactionsICs}
\end{equation}
In analogy to the porous medium equation \eqref{eq:pme}, a crucial question is whether there are (weak) solutions of this model that feature a propagating interface.  In the context of \eqref{eq:pdeWithReactions}, an interface is a point at which the resident population $\tilde u$ is at its carrying capacity while the invading population $v$ is at zero.  A schematic of this setup is shown in Figure~\ref{fig:schematic}(b).

The existence or otherwise of such an interface has not been conclusively settled by previous studies.  On the one hand, there are a number of examples of numerical simulations of PDEs of the form \eqref{eq:pdeWithReactions} with initial conditions \eqref{eq:pdeWithReactionsICs} that do appear to have a propagating interface, often evolving to travelling wave solutions (depending on the form of the reaction terms $F$ and $G$) \cite{colson2021travelling, el2021travelling, mascia2024numerical}.  However, in \cite{gallay2022propagation} it was shown that a one-parameter family of travelling waves to their model must exhibit algebraic decay of $V$ in front of the wave, rather than an interface at a finite location.  The front of the wave becomes steeper in the limit that the speed of the wave is taken to zero.  In addition, numerical estimates of the wave speed from PDE simulations made in \cite{mascia2024numerical} indicated minor but non-negligible dependence on the spatial resolution of the numerical scheme, raising concerns as to whether the results of the PDE simulations had in fact numerically converged.

In this article we consider in depth an idealised version of systems of the form \eqref{eq:pdeWithReactions} that includes the coupled effects of degenerate cross-diffusion and degradation:
\begin{equation}
\pfrac{u}{t} = v, \qquad \pfrac{v}{t} = \pfrac{}{x}\left(u\pfrac{v}{x}\right).
\label{eq:pde}
\end{equation}
Here we have defined the resident population deficit $u = 1-\tilde u$, so that $u$ may be considered the (time and space-dependent) diffusivity itself, and $u=0$ represents the resident population at capacity.  We are particularly interested in whether we can find solutions for which the diffusivity is initially zero, and the invading population is compactly supported, such as in \eqref{eq:pdeWithReactionsICs}.  However, an initially compactly supported $v$ will (in stark contrast to the porous media equation \eqref{eq:pme}) prove to be problematic.  We thus consider the initial condition
\begin{equation}
u(x,0) = 0, \qquad v(x,0) = \begin{cases} 1 & x < x^*_0 \\ \epsilon & x > x^*_0 \end{cases},
\label{eq:pdeIC}
\end{equation}
where the initial height in front of the discontinuity $\epsilon \in (0,1)$.  To examine the possibility of a solution with an interface, we will then focus on the limiting behaviour as this leading height, $\epsilon$, becomes small.

Our main finding will be that the limit $\epsilon \to 0$ is singular.  While for small $\epsilon$ an apparent interface (a point at which $v$ is small and the slope of $v$ is large in $\epsilon$) will be present, the speed of propagation of this point vanishes, and the slope of the profile becomes infinite, as $\epsilon \to 0$.  Our conclusion is thus that the system \eqref{eq:pde} does not support solutions that have a moving interface, in contrast to systems such as the porous medium equation \eqref{eq:pme}.  
However, this limit turns out to be approached logarithmically slowly in $\epsilon$ (to be precise, $O([\log(1/\epsilon)]^{-1/2})$), which explains why numerical simulations of such models may appear to converge, as they have numerically regularised the problem at $O(\epsilon)$, while observed quantities such as speeds of propagation superficially appear to be order one.

We start in Section~\ref{sec:numerics} by showing an example numerical simulation of \eqref{eq:pde}, and demonstrating how solutions for any given $\epsilon > 0$ can be found by solving a similarity solution rather than the full PDE.  The numerical calculation of these similarity solutions shows a number of interesting properties, including the possibility of persistent nonexistence of spatial derivatives of $u$ (and higher spatial derivatives of $v$) over time, as well as showing the first indication of logarithmically slow behaviour of solution parameters (such as the speed of propagation) in the $\epsilon \to 0$ limit.  In Section~\ref{sec:asymptotics} we perform a matched asymptotic analysis of the similarity solution for small $\epsilon$, which confirms this logarithmic behaviour.  We discuss the more general implications of these results and future extensions in Section~\ref{sec:discussion}.

\begin{figure}
\centering
\includegraphics[width=\textwidth]{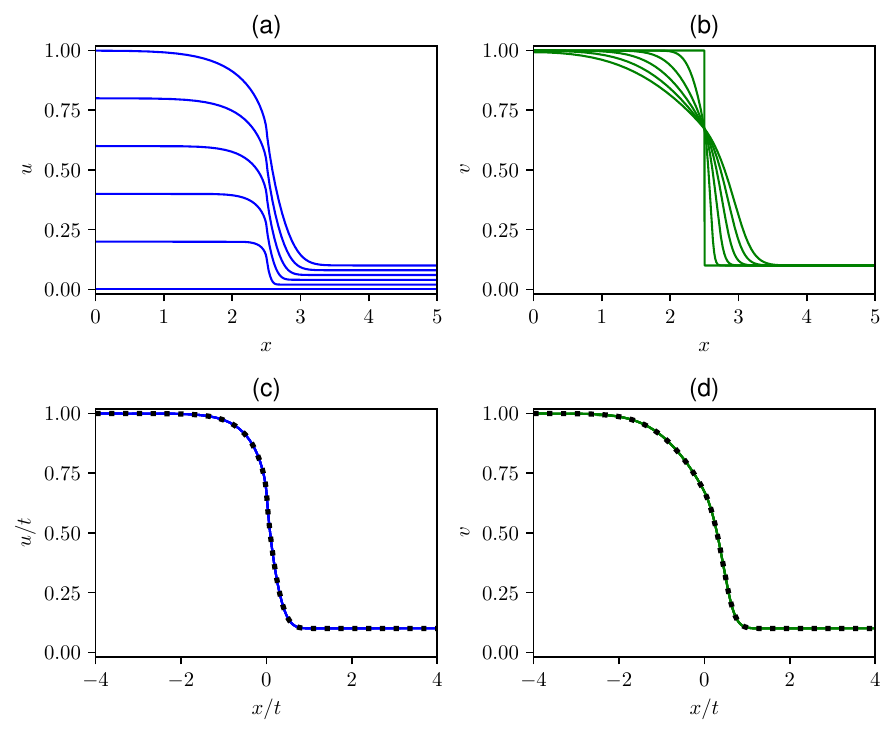}
\caption{Numerical solutions of the PDE system \eqref{eq:pde} for (a) $u(x,t)$ and (b) $v(x,t)$, at evenly spaced times $t$ from $0$ to $1$.  The initial condition is \eqref{eq:pdeIC} for $\epsilon = 0.1$.  In (c), (d), the PDE simulations at the times plotted in (a,b) (except for the initial condition $t=0$) are scaled onto similarity variables $\eta = (x-x^*_0)/t$, $f = u/t$, $g = v$, showing exact collapse to a similarity solution.  Black dotted curves are the corresponding numerical solution of the ODE system \eqref{eq:simode}.}
\label{fig:solutionProfiles}
\end{figure}

\section{Numerical solutions and self-similarity}
\label{sec:numerics}

To demonstrate the generic behaviour of the PDE system \eqref{eq:pde} with initial condition \eqref{eq:pdeIC}, we simulate the system numerically.  This is performed on a finite domain $0 < x < 5$, with no-flux conditions on $v$ at the ends (no boundary conditions are required for $u$ as it has no spatial derivatives).  The initial discontinuity is taken at the centre of the domain at $x^*_0 = 2.5$.  The PDE is discretised in space using a finite volume method implemented in MATLAB, and advanced in time using MATLAB's \texttt{ode15s} implicit time stepping algorithm.  The details of this method are quite standard and very similar to those used in \cite{chen2026wavespeed}, so we do not go into further detail here.

We show the results of such a simulation for $\epsilon=0.1$ in Figure~\ref{fig:solutionProfiles}(a,b), for $u$ and $v$, respectively.  Solution profiles are shown for equally spaced times from $t=0$ to $t=1$.  This solution demonstrates the coupled behaviour of the two variables, with $u$ increasing more rapidly to the left of the initial discontinuity, and $v$ diffusing in both directions, although more rapidly to the left of the initial discontinuity and more slowly to the right.  This behaviour is as expected given that $u$ is the diffusivity of $v$.

\subsection{Formulation of system for similarity solution}

While the system can be simulated as above for any $\epsilon$, a chief advantage of considering the simplified model \eqref{eq:pde} with initial condition \eqref{eq:pdeIC} is that it permits a similarity solution.  This similarity solution may be calculated numerically by solving a single boundary value problem with much greater precision than the original PDE, which will be particularly valuable when looking at small values of the height $\epsilon$ in front of the initial discontinuity.  We also note that, while exact for the simplified model \eqref{eq:pde}, this similarity solution also represents the early-time behaviour of more complex cross-diffusion models of the form \eqref{eq:pdeWithReactions} evolving from an initial condition \eqref{eq:pdeIC}, and so captures essential behaviour of this class of models as well (see Section \ref{sec:discussion} for further discussion).

To construct similarity solutions we first define the similarity variable $\eta$, and similarity profiles $f$ and $g$, by
\begin{equation}
u(x,t) = tf(\eta), \qquad v(x,t) = g(\eta), \qquad \eta = \frac{x-x^*_0}{t}.
\label{eq:similarityAnsatz}
\end{equation}
The factors $t$ in \eqref{eq:similarityAnsatz} result from requiring the terms in \eqref{eq:pde} have the same dependence on $t$, while the scale of $v$ cannot change as it must match to the fixed values behind and ahead of the initial discontinuity in the large $\eta$ limit.  From \eqref{eq:pde} the profiles $f$ and $g$ then satisfy the ODE system
\begin{equation}
f - \eta \ddfrac{f}{\eta} = g, \qquad -\eta \ddfrac{g}{\eta} = \ddfrac{}{\eta}\left(f\ddfrac{g}{\eta}\right), \qquad -\infty < \eta < \infty
\label{eq:simode}
\end{equation}
with far-field conditions consistent with the initial condition \eqref{eq:pdeIC}:
\begin{equation}
f, g \to \begin{cases} 1 & \eta \to -\infty \\ \epsilon & \eta \to \infty. \end{cases}
\label{eq:odeBCs}
\end{equation}
As $t\to 0^+$, the profiles contract around the point $x=x^*_0$, thus a similarity solution satisfying these conditions will approach the correct initial condition \eqref{eq:pdeIC}.

Looking closer at the far-field conditions \eqref{eq:odeBCs}, if we expand $f$ and $g$ close to a constant in each far field $|\eta|\to\infty$, we find that \eqref{eq:simode} supports the behaviour
\begin{equation}
f \sim C_0 + C_1\eta + C_2(\eta) \exp\left(-\frac{\eta^2}{2C_0}\right), \qquad |\eta| \to\infty.
\label{eq:farFieldBehaviour}
\end{equation}
for $C_0$, $C_1$ constants, and the WKB prefactor $C_2$ a function of $\eta$ (whose exact form is unimportant).  In order to satisfy the far-field conditions we require $C_0 = 1$ for $\eta \to -\infty$ and $C_0 = \epsilon$ for $\eta\to \infty$, while the algebraically growing mode $C_1$ must be zero in both cases.  The remaining WKB mode decays exponentially in each case.  Thus each far-field condition effectively applies two constraints on the system, resulting in four boundary conditions for the third-order system \eqref{eq:simode}.  

In order to find solutions then, an additional degree of freedom is required.  This is provided by the fact that the system \eqref{eq:simode} is singular at the point at $\eta = 0$.  Expanding $f$ and $g$ near $\eta=0$ we find that $f$ may be nondifferentiable; in particular,
\begin{equation}
f = a_0 + a_1\eta\log|\eta| + a_2\eta + \ldots, \qquad g = a_0 - a_1\eta + \ldots,
\label{eq:expansionNearZero}
\end{equation}
where $a_0$, $a_1$, $a_2$ are arbitrary, and all higher coefficients may be determined in terms of these three constants.
For $g$ to be differentiable across the singular point, $a_0$ and $a_1$ must take the same value on either side of $\eta=0$.  However, $a_2$ may take different values on either side, say
\[
a_2 = \begin{cases} a_2^- & \eta < 0 \\ a_2^+ & \eta > 0. \end{cases}
\]
This extra degree of freedom allows us to satisfy both far field conditions.  A consequence of the singular point is that even for time greater than zero, $u$ will be nondifferentiable at the point $x_0$ at which the initial discontinuity existed.  The PDE system \eqref{eq:pde} permits this because there is no diffusion of $u$.

\subsection{Numerical calculation of similarity solutions}

We now describe the method of numerically calculating solutions of \eqref{eq:simode}.  To achieve high accuracy, and observe the change in profiles as $\epsilon$ is varied, we use the numerical continuation package AUTO-07p \cite{doedel2007auto}.  To obtain a solution with the appropriate boundary conditions \eqref{eq:farFieldBehaviour}, we will apply a homotopy continuation method, in which an initial trial solution is computed with explicit error terms in the boundary conditions, which are then each reduced to zero in turn via numerical continuation.  This process has been applied to the computation of similarity solutions previously with much success~\cite{dallaston2017self,tseluiko2013homotopy}.

To correctly treat the singular point at $\eta=0$, we formulate the system as two sets of ODEs, one starting at $\eta = -\Delta\eta$, just to the left of the singular point, and extending to $\eta = -L$ in the negative direction, and one starting from $\Delta\eta$ and extending to $\eta = R$ in the positive direction.  Here $\Delta\eta \ll1$ and $R,L \gg 1$; in our calculations we take $\Delta\eta = 10^{-4}$, $R=8$, and $L=5$; the exponentially decaying term in \eqref{eq:expansionNearZero} means that the far field condition is approached rapidly in each direction, so $L$ and $R$ do not have to be very large.  
Defining the scaled flux variable $h(\eta) = u(\eta) v'(\eta)$, we then have the system
\begin{subequations}
\begin{align}
\ddfrac{f^\pm}{\eta} &= \frac{f^\pm - g^\pm}{\eta} \\
\ddfrac{g^\pm}{\eta} &= \frac{h^\pm}{f^\pm} \\
\ddfrac{h^\pm}{\eta} &= -\eta \frac{h^\pm}{f^\pm},
\end{align}
\label{eq:odeSystem}%
\end{subequations}
where $\pm$ superscripts refer to the negative and positive $\eta$-directions, respectively.

Close to the singular point we use the expansion \eqref{eq:expansionNearZero} as boundary conditions:
\begin{equation}
\left.
\begin{array}{ll}
f^\pm = a_0 \pm a_1 \Delta\eta \log(\Delta\eta) \pm a_2^\pm \Delta\eta \\
g^\pm = a_0 \pm a_1\Delta\eta\\
h^\pm = a_0a_1,
\end{array}
\right\},
\quad \eta = \pm\Delta\eta.
\end{equation}
Here $a_0$, $a_1$, and $a_2^\pm$ are parameters that will be varied in order to satisfy the far field conditions.  In the far field we apply
\begin{align}
f^+ - g^+ &= \mathcal E^+, \qquad f^+ = C_0^+, & \eta = R \\
f^- - g^- &= \mathcal E^-, \qquad f^- = C_0^-, & \eta = -L.
\end{align}  
Here $\mathcal E^{\pm}$ represent error terms from the desired boundary conditions \eqref{eq:odeBCs} (in which $f$ and $g$ are equal), while $C_0^\pm$ represents each far field value.  The aim of the homotopy continuation procedure will be to reduce each of $\mathcal E^\pm$ to zero, set $C_0^- = 1$, and $C_0^+ = \epsilon$ for any given $\epsilon$.

We proceed as follows.  We start with the exact linear solution to \eqref{eq:simode}, $f = \eta/2$, $g = 1$, with domain size parameters $R = L = 1$ (the parameter $\Delta\eta =  10^{-4}$ is held constant throughout the continuation procedure).  This initial solution then satisfies the imposed boundary conditions given the parameter values:
\[
a_0 = 1, \quad a_1 = 0, \quad a_2^\pm = \mp \frac{1}{2}, \quad C_0^- = \frac{3}{2}, \qquad C_0^+ = \frac{1}{2}, \qquad \mathcal E^\pm = \mp \frac{1}{2}.
\]
We then apply numerical continuation, firstly reducing the left-hand error $\mathcal E^-$ to zero and the left-hand value $C_0^-$ to one, while  letting $a_0, a_2^-$, and $\mathcal E^+$ vary.  We then reduce the right hand error $\mathcal E^+$ to zero and $C_0^+$ to $\epsilon=0.1$, letting each of $a_0, a_1, a_2^\pm$ vary.  Finally, we increase the domain size $L$, $R$ until the domain size has negligible effect on the solution.  At the end of this process, we have a similarity solution for the given $\epsilon$.  We may then numerically continue in $\epsilon$ to find the solution for any such $\epsilon$.

In Figure~\ref{fig:solutionProfiles}(c,d) we plot the resulting similarity solution for $\epsilon = 0.1$.  Comparing this similarity solution to the PDE solution profiles collapsed onto the self-similar coordinates, we see the profiles agree exactly, thus confirming the similarity solution is an exact solution of the original PDE for the given initial condition.  For the remainder of the study, we thus focus purely on the calculation of the similarity solution itself.

\subsection{Numerical results for small $\epsilon$}
\label{sec:numericalResultsForSmallEpsilon}
\begin{figure}
\centering
\includegraphics[width=\textwidth]{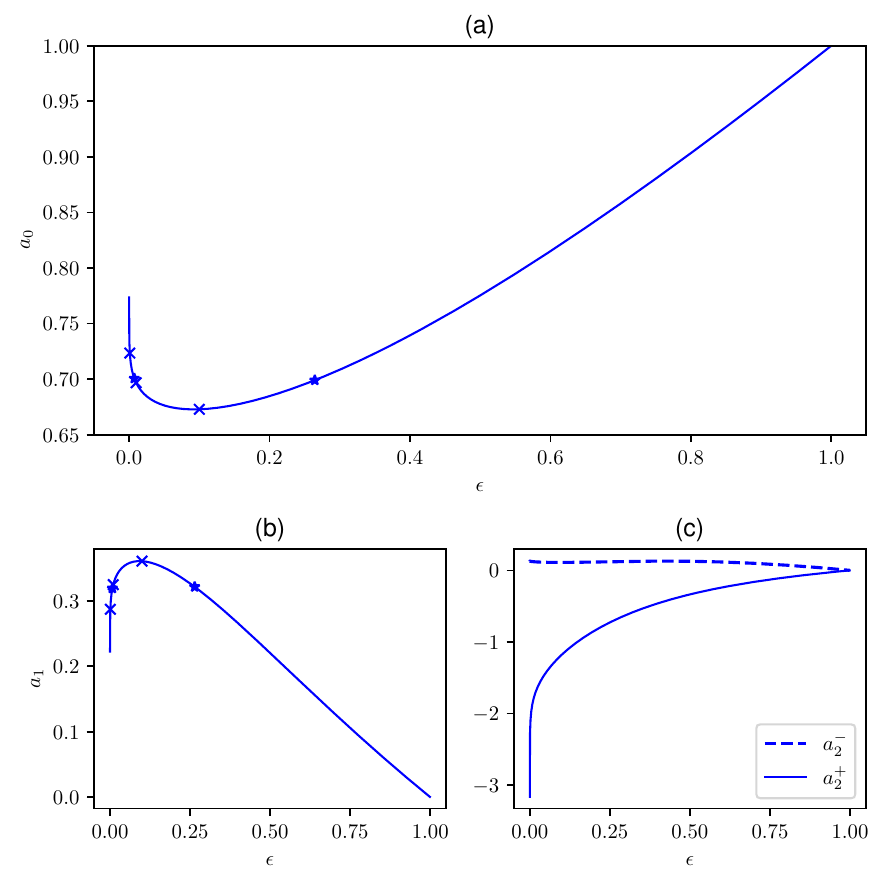}
\caption{Coefficients from the expansion \eqref{eq:expansionNearZero} of similarity solutions  \eqref{eq:simode} over front height $\epsilon$:  (a) the value $a_0$, (b) $a_1$, and (c) $a_2^\pm$. The parameters $a_0$ and $a_1$ show a turning point, below which $a_0$ and $a_1$ appear to slowly return to $1$ and $0$, respectively.  The value of $a_2^+$ increases in magnitude (slowly) as $\epsilon \to 0$.  The crosses marked on the curves in (a,b) correspond to the values $\epsilon = 10^{-1}, 10^{-2}, 10^{-3}$ shown in Figure~\ref{fig:exampleSolutions}(a,b), while the asterisks correspond to the values $\epsilon = 0.268, 0.0078$ shown in Figure~\ref{fig:exampleSolutions}(c,d).}
\label{fig:epsilonCurves}
\end{figure}

Now having a method of calculating similarity solutions, we can readily examine how solutions depend on the height $\epsilon$.  In Figure~\ref{fig:epsilonCurves} we show how the constants $a_0$ and $a_1$ depend on $\epsilon$ over the range $\epsilon \in (0,1)$.  These constants may be considered characteristic parameters describing a solution to the system \eqref{eq:pde}, given
\[
u(x^*_0,t) = v(x^*_0,t) = a_0, \qquad \pfrac{v}{x}(x^*_0,t) = -\frac{a_1}{t},
\]
where $x^*_0$ is the location of the initial discontinuity.

\begin{figure}
\centering
\includegraphics[width=\textwidth]{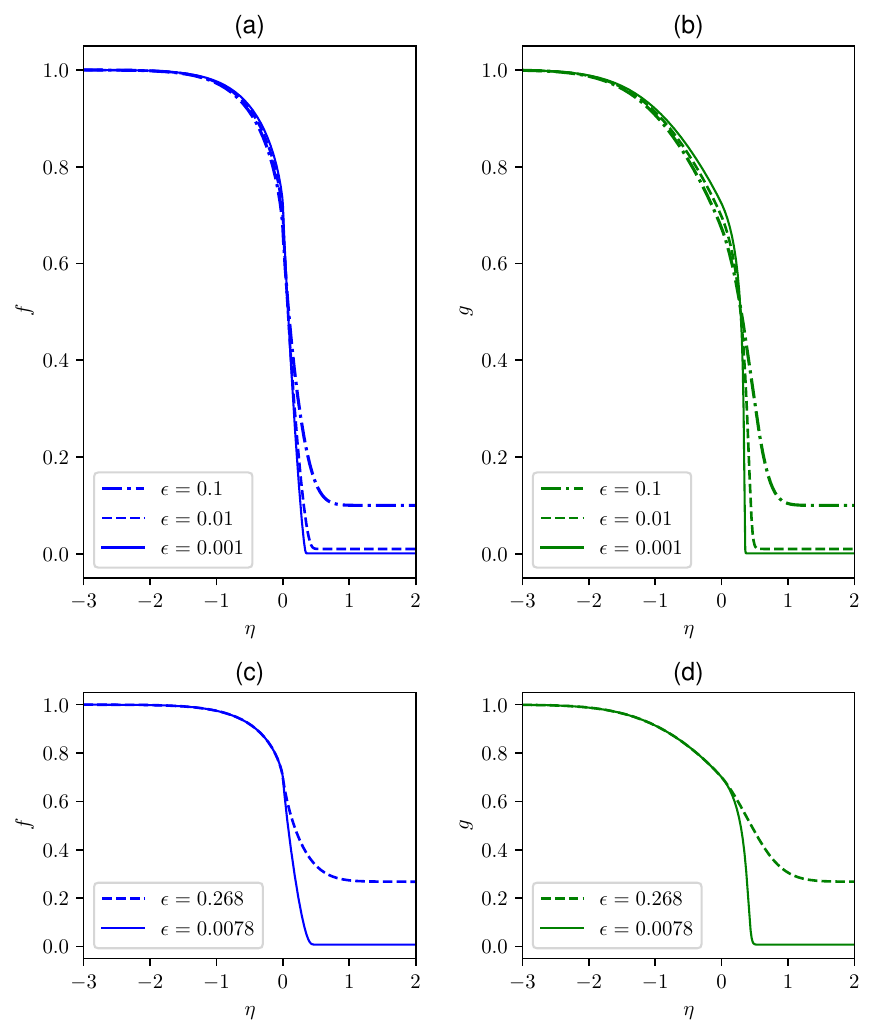}
\caption{(a,b) Similarity profiles $f$ and $g$, respectively, as $\epsilon$ becomes small.  An apparent interface at $\eta=\eta^*$ is seen, where the magnitude of $f$ becomes of order $\epsilon$ and rapidly flattens out.  In (c,d), two solutions for different values may have the same value of $a_0$, and thus the same profile for $\eta < 0$.}
\label{fig:exampleSolutions}
\end{figure}

As $\epsilon$ decreases, at first $a_0$ decreases and $a_1$ increases.  However, there is a turning point at around $\epsilon = 0.1$ below which $a_0$ begins to increase and $a_1$ to decrease.   Focusing on small $\epsilon$, in Figure~\ref{fig:exampleSolutions}(a,b) we plot solution profiles for $\epsilon = 10^{-1}, 10^{-2}, 10^{-3}$.
For small $\epsilon$, the similarity solutions are qualitatively different on either side of the singular point $\eta=0$.  In particular, the solution profiles are quite flat to the left of the singular point, but to the right, a steep front forms.  Close to the bottom of the front the solution reaches a point where the slopes of both $f$ and $g$ rapidly flatten as $f$ and $g$ approach $\epsilon$.  On the other hand, the values of $a_2^\pm$ differ (see Figure~\ref{fig:epsilonCurves}(c)), in particular with $a_2^+ < 0$ increasing in magnitude as $\epsilon$ becomes small.  This contrast indicates the magnitude of derivatives of a solution to the left and the right of the singular point becomes markedly different.

The most notable feature for solutions when $\epsilon$ is small is the formation of what could be considered an apparent interface at a location $\eta^* > 0$.  This is not a true interface, as for $\epsilon$ positive $v\to \epsilon$ never vanishes, but the slope of $g$ appears almost vertical near this point, and the change in slope of $f$ and $g$ occurs very rapidly.  In a solution to the PDE \eqref{eq:pde}, this apparent interface will occur at $x^*(t) = x_0^* + \eta^*t$, thus the value of $\eta^*$ in the similarity solution represents the speed of propagation of such a point.

\begin{figure}
\centering
\includegraphics[width=\textwidth]{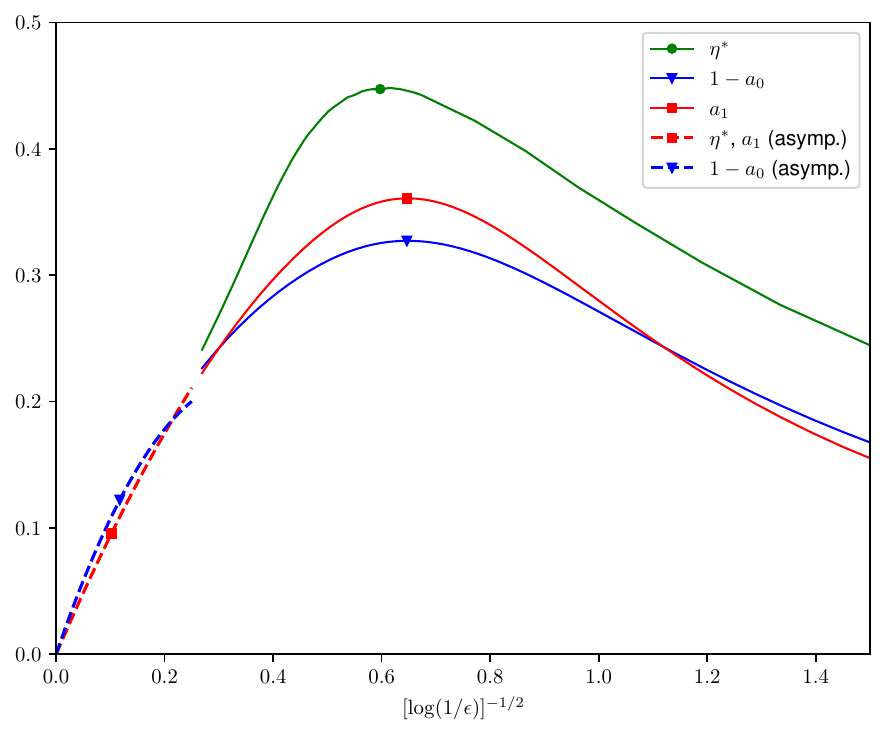}
\caption{Apparent interface location $\eta^*$, as well as $1-a_0$ (where $a_0 = f(0) = g(0)$) and slope $a_1 = -g'(0)$, in the limit that $\epsilon$ becomes small.  Numerical results indicate that these solution parameters are all going to zero, although only logarithmically fast, namely, of order $[\log(1/\epsilon)]^{-1/2}$ as $\epsilon \to 0$.  Also plotted are the two-term asymptotic approximations \eqref{eq:finalSeries} for each parameter (note these approximations are valid in the limit $\epsilon \to 0$).}
\label{fig:epsilonCurvesEtaS}
\end{figure}

To characterise the location of $\eta^*$ over $\epsilon$, we define it as the point at which the slope of $g$ becomes smallest (that is, most negative).  In Figure~\ref{fig:epsilonCurvesEtaS} we plot this quantity, along with the difference between $a_0$ and its asymptotic value $1$, and the slope $a_1$.  As each of these quantities is seen to tend to zero as $\epsilon\to 0$, but slowly, we plot these against the logarithmic variable $[\log(1/\epsilon)]^{-1/2}$.  Here we see initial signs that these variables are of this logarithmic order as $\epsilon \to 0$, which will motivate the asymptotic analysis in the next section (for comparison we also plot the results of the asymptotic analysis we will undertake).  Note that in this plot we take $\epsilon$ down to $10^{-6}$, which is at about the limit of the numerical solution procedure.  For this value, $[\log(1/\epsilon)]^{-1/2} \approx 0.27$ is still relatively large, which increases the difficulty of deducing the behaviour in this limit.

Finally we note some interesting implications of the nonmonotonic nature of the constants $a_0$ and $a_1$ seen in Figure~\ref{fig:epsilonCurves}.  One implication, which can be seen in the small $\epsilon$ solution in Figure~\ref{fig:exampleSolutions}(a, b), is that a smaller value of $\epsilon$ to the right of the singular point may result in a larger value of $g$ (and thus $v$) to the left of the singular point.  This means that given initial conditions for which the diffusivity $u$ vanishes, the problem \eqref{eq:pde} does not necessarily satisfy a comparison principle; one may have two solutions $v_1$ and $v_2$, where $v_1 \leq v_2$ at $t=0$, but $v_1$ does not remain bounded by $v_2$ for $t > 0$.

The profile to the left of $\eta=0$ is uniquely characterised by $a_0$, as for a given $a_0$, both $a_1$ and $a_2^-$ are determined by requiring the $\eta\to-\infty$ condition to be met.  A second implication of nonmonotonic $a_0$ is thus that two different values of $\epsilon$ can result in two similarity solutions that are identical for $\eta <0$, while being distinct for $\eta > 0$.  In Figure~\ref{fig:exampleSolutions}(c, d) we plot the similarity solutions for $\epsilon \approx 0.268$ and $\epsilon \approx 0.0078$, both corresponding to a value of $a_0 = 0.7$.  Here we see the exact agreement for $\eta < 0$.  Given the similarity solutions represent exact solutions of the pde, this implies that two initial conditions to \eqref{eq:pde} which are the same on one side of the singular point $x = x^*_0$ could remain the same for $t>0$, even if this is not the case on the other side of the singular point.

\section{Asymptotic analysis for small $\epsilon$}
\label{sec:asymptotics}

As indicated by the numerical solutions, solution characteristics, including the apparent interface location $\eta^*$, appear to very slowly (logarithmically) tend to zero as the front height $\epsilon$ vanishes.  In this section, we establish this result systematically using a matched asymptotic analysis.

\begin{figure}
\centering
\includegraphics[]{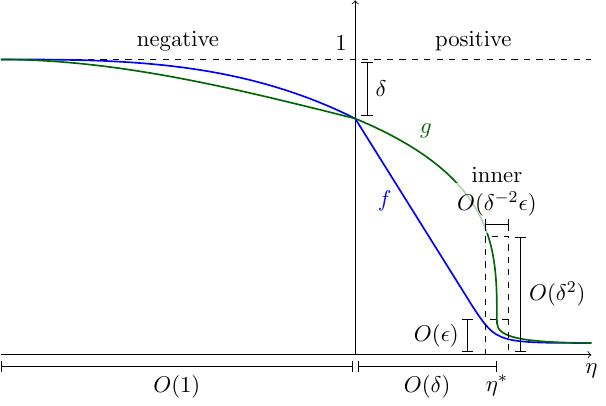}
\caption{Schematic of the asymptotic structure of similarity solutions to \eqref{eq:simode} as described in Section \ref{sec:asymptotics}.  The solution has three main regions: the negative region $\eta < 0$, the positive region $0 < \eta < \eta^*$ (where $\eta^*$ is the apparent front location), and an inner region close to $\eta^*$.  Matching in the inner region reveals that $\delta = O(\log(1/\epsilon)^{-1/2})$.}
\label{fig:asymptoticSchematic}
\end{figure}

In Figure~\ref{fig:asymptoticSchematic} we sketch an overview of the asymptotic structure that we uncover in this section.  The structure of the solution profiles in the small-$\epsilon$ limit are starkly different for $\eta < 0$, left of the singular point, and $\eta > 0$ to the right. For $\eta < 0$ both $f$ and $g$ become close to flat, whereas for $\eta > 0$, the solution profiles are compressed in a region $0 < \eta < \eta^*$, which becomes logarithmically small in the limit $\epsilon \to 0$.  Finally, there is an inner boundary layer problem close to $\eta^*$ that must be resolved, and matched to the outer region, in order to determine the relationship between each solution parameter.

We start by defining a small parameter $\delta$, by
\begin{equation}
\delta = 1-a_0.
\label{eq:deltaDefinition}
\end{equation}
As seen in Figure~\ref{fig:epsilonCurvesEtaS}, we expect $\delta$ to be logarithmically small as $\epsilon$ vanishes; the goal of the asymptotic analysis will be to determine this relationship precisely.  We will calculate as many terms as required to find two terms in the series approximation of $\delta(\epsilon)$.  Higher terms quickly become infeasible due to the complexity of the integrations required.

\subsection{Negative $\eta$}

When $\eta < 0$, the profiles $f$ and $g$ are both constrained in magnitude by the far-field condition $f, g \to 1$ as $\eta\to -\infty$.  Given the definition \eqref{eq:deltaDefinition}, we assume the asymptotic expansion
\begin{equation}
f(\eta) = 1 + \delta f_1(\eta) + \delta^2 f_2(\eta) + \ldots, \qquad g(\eta) = 1 + \delta g_1(\eta) + \delta^2 g_2(\eta) + \ldots
\end{equation}
Substituting into the ODE system \eqref{eq:simode} and taking successive terms in $\delta$ we have
\begin{equation}
\ddfrac{^2g_1}{\eta^2} + \eta \ddfrac{g_1}{\eta} = 0, \qquad 
\ddfrac{^2g_2}{\eta^2} + \eta \ddfrac{g_2}{\eta} = -\ddfrac{}{\eta}\left(f_1\ddfrac{g_1}{\eta}\right)
\end{equation}
and so on, with each $f_j$ found from $g_j$ by integrating:
\[
f_j = -\eta \int \frac{g_j}{\eta^2} \, \mathrm d\eta, \qquad j = 1, 2, \ldots.
\]
To satisfy the far-field condition we require $f_j, g_j\to 0$ as $\eta \to -\infty$, while from the definition of $\delta$:
\[
f_j(0), g_j(0) = \begin{cases} -1 & j = 1 \\ 0 & j = 2,3,\ldots \end{cases}
\]
We thus have
\begin{subequations}
\begin{equation}
g_1 = -\left(1 + \erf\left(\frac{\eta}{\sqrt 2}\right)\right), \qquad f_1 = \frac{\eta}{\sqrt{2\pi}}\Ei\left(-\frac{\eta^2}2\right) - 1 - \erf\left(\frac{\eta}{\sqrt{2}}\right),
\label{eq:negativeSolutionFirst}
\end{equation}
where $\erf$ is the error function, and $\Ei$ is the exponential integral, defined for negative arguments by
\[
\Ei(s) = \int_{-\infty}^s \frac{\e^{s'}}{s'} \, \mathrm ds'.
\]
The next term $g_2$ is complicated but possible to find explicity by integrating:\begin{multline}
g_2 = \left(1+\erf\left(\frac{\eta}{\sqrt{2}}\right)\right) \left(-\frac{e^{-\frac{\eta^2}{2}} \eta}{\sqrt{2 \pi }}-\frac{1}{3}+\frac{1}{3\pi }-\frac{\log (2)}{3 \pi }\right) -\frac{1}{3}\erf\left(\frac{\eta}{\sqrt{2}}\right)^2 \\
+\frac{e^{-\frac{\eta^2}{2}} (\eta^2-1)}{3 \pi }\Ei\left(-\frac{\eta^2}{2}\right) +\frac{\Ei\left(-\eta^2\right)}{3 \pi }-\frac{e^{-\eta^2}}{3 \pi }+\frac{1}{3}.
\end{multline}
\label{eq:negativeExpansionSolution}%
\end{subequations}
The derivative of $g_2$ at $\eta=0$ is needed to give the correction to the slope $a_1$ that matches with that for $\eta>0$.  Differentiating the above we find that
\begin{equation}
a_1 = a_{11}\delta + a_{12}\delta^2 + \ldots,
\end{equation}
where
\begin{equation}
a_{11} = -\left.\ddfrac{g_1}{\eta}\right|_{\eta\to0^-} = \sqrt{\frac{2}{\pi}}, \qquad a_{12} = -\left.\ddfrac{g_2}{\eta}\right|_{\eta\to0^-} = \frac{5\pi - 2 - 2\log(2)}{3\sqrt{2}\pi^{3/2}}.
\label{eq:a1Terms}
\end{equation}

\subsection{Positive $\eta$ away from apparent interface}

We now consider the interval $0 < \eta < \eta^*$, where $\eta^*$ is the apparent interface location.  Since $\eta^*$ is seen to decrease at the same logarithmic rate as $\delta$, we introduce the scaled independent variable $\xi$ by $\eta = \delta\xi$, so that \eqref{eq:simode} becomes
\begin{equation}
f - \xi\ddfrac{f}{\xi} = g, \qquad \ddfrac{}{\xi}\left(f\ddfrac{g}{\xi}\right) = -\delta^2\xi\ddfrac{g}{\xi}.
\label{eq:odeInXi}
\end{equation}
We take the apparent interface location $\eta^* = \delta\xi^*$, where $\xi^*$ has a series expansion in $\delta$:
\begin{equation*}
\xi^* = \xi_0^* + \delta \xi_1^* + \ldots
\end{equation*}
Expanding the similarity profiles $f$ and $g$ in this region as
\begin{equation*}
f = f_0(\xi) + \delta f_1(\xi) + \delta^2 f_2(\xi) + \ldots, \quad g = g_0(\xi) + \delta g_1(\xi) + \delta^2 g_2(\xi) + \delta^3 g_3(\xi) + \ldots,
\end{equation*}
the first terms must satisfy the differential equations
\begin{subequations}
\begin{align}
&f_0 - \xi \ddfrac{f_0}{\xi} = g_0, \qquad \ddfrac{}{\xi}\left(f_0\ddfrac{g_0}{\xi}\right) = 0 \\
&f_1 - \xi \ddfrac{f_1}{\xi} = g_1, \qquad \ddfrac{}{\xi}\left(f_0\ddfrac{g_1}{\xi} + f_1\ddfrac{g_0}{\xi}\right) = 0.
\end{align}
At $\xi=0$, the value and derivative of $g$ must match with those for $\eta \to 0^-$, giving the boundary conditions
\begin{equation}
g_0 = 1, \quad g_1= -1, \quad \ddfrac{g_0}{\xi}= 0, \quad \ddfrac{g_1}{\xi} = 0, \qquad \xi = 0.
\end{equation}
Another boundary condition comes from the apparent interface at $\xi^*$.  We expect $f$ to be exponentially small in $\delta$ close to this point, so expanding $\xi^*$ and $f$, each term in $\delta$ must vanish:
\begin{equation}
f_0(\xi_0^*) =0, \qquad f_1(\xi_0^*) + \left.\ddfrac{f_0}{\xi}\right|_{\xi=\xi^*}\xi_1^* = 0.
\end{equation}
\label{eq:outerProblemFirst}
\end{subequations}
This set of equations and boundary conditions implies that $f_0$ and $f_1$ are linear, while $g_0$ and $g_1$ are constant:
\begin{equation}
g_0 = 1, \quad g_1 = -1, \quad f_0 = 1 - \frac{\xi}{\xi_0^*}, \quad f_1 = -1 + \frac{\xi_0^* + \xi_1^*}{\xi_0^{*2}}\xi.
\label{eq:outerSolutionFirst}
\end{equation}

We can continue to higher terms.  We will calculate the next corrections $g_2$ and $g_3$, which satisfy
\begin{subequations}
\begin{equation}
\ddfrac{}{\xi}\left(f_0\ddfrac{g_2}{\xi}\right) = 0, \qquad
\ddfrac{}{\xi}\left(f_0\ddfrac{g_3}{\xi}\right) = -\ddfrac{}{\xi}\left(f_1\ddfrac{g_2}{\xi}\right)
\end{equation}
with boundary conditions
\begin{equation}
g_2 = g_3 = 0, \quad \ddfrac{g_2}{\xi} = -a_{11}, \quad \ddfrac{g_3}{\xi} = -a_{12}, \qquad \xi=0.
\end{equation}
\label{eq:outerProblemNext}%
\end{subequations}
Note that the higher corrections to $f$ are not required to calculate $g$ to this order.  Integrating and applying the conditions at zero, we find
\begin{subequations}
\begin{align}
g_2 &= a_{11} \xi_0^* \log\left(1- \frac{\xi}{\xi_0^*}\right) \\
 g_3 &= a_{11}\frac{\xi_1^*\xi}{\xi_0^*-\xi} + \left(a_{12}\xi_0^* + a_{11}\xi_1^*\right)\log\left(1-\frac{\xi}{\xi_0^*}\right).
\end{align}
\label{eq:outerSolutionNext}%
\end{subequations}

\subsection{Inner problem near apparent interface}

The expansions \eqref{eq:outerSolutionFirst}, \eqref{eq:outerSolutionNext} are singular in the limit $\xi \to \xi^*$.  There is thus a region close to $\xi=\xi^*$ where the solution will match onto an inner solution, in which the far field condition $f \to \epsilon$ is attained.
We thus define inner variables $z$, $\hat f$, and $\hat g$, by 
\[
\xi = \xi^* + \delta^{-2}\epsilon z, \qquad f = \epsilon \hat f(z), \qquad g = \delta^2 \hat g(z).
\]
Here, the scale of $f$ is determined so that $\hat f \to 1$ as $z\to\infty$, while the scales for $z$ and $g$ are determined by finding the appropriate balance of terms in \eqref{eq:odeInXi}.  Under this change of variables, the system \eqref{eq:odeInXi} is, to leading order,
\begin{equation}
-\xi^* \ddfrac{\hat f}{ z} = \hat g, \qquad -\xi^*\ddfrac{\hat g}{z} = \ddfrac{}{z}\left(\hat f \ddfrac{\hat g}{z}\right).
\label{eq:innerProblem1}
\end{equation}
Note that  a result of solving the system in a small region near $\xi^*$ is that the system is autonomous.  Solutions of \eqref{eq:innerProblem1} are effectively travelling wave solutions of \eqref{eq:pde}; this inner problem is thus relevant close to an apparent interface not just for the similarity solution but in more general contexts (see further discussion in Section \ref{sec:discussion}).

To match to the far field condition, we require $\hat f \to 1$, $\hat f' \to 0$ as $z\to \infty$.  Integrating \eqref{eq:innerProblem1} twice and applying these conditions, we find
\begin{equation}
\ddfrac{\hat f}{z} = -\xi^*\log(\hat f),
\label{eq:innerProblem2}
\end{equation}
thus we have an implicit solution as follows:
\begin{equation}
z = z_0 - \frac{1}{\xi^*}\li(\hat f) = z_0 - \frac{1}{\xi^*}\Ei(\log(\hat f)), \qquad \hat g = \xi^{*2}\log(\hat f)
\label{eq:innerSolution}
\end{equation}
where $z_0$ is an arbitrary translational constant, and $\li$ is the logarithmic integral.  Note that $\hat g$ vanishes in the limit $z \to\infty$, which is consistent with the definition of $\hat g$ given $\delta^2 \gg\epsilon$; an additional inner problem would be required to match $g$ onto its far-field condition $g \to \epsilon$.

What remains is to find the asymptotic relationship between $\delta$ and $\epsilon$, as well as the asymptotic formula for the apparent interface location $\xi^*$.  To do so we need to match the $z\to -\infty$ expansion of the inner problem \eqref{eq:innerSolution} to the $\xi\to\xi^*$ expansion of the outer problem \eqref{eq:outerSolutionFirst}, \eqref{eq:outerSolutionNext}.

While the asymptotic expansion for the inverse logarithmic integral could be applied directly to \eqref{eq:innerSolution}, a simpler way to carry out the matching is to write the relationship between $\hat g$ and $\hat f$ in \eqref{eq:innerSolution} in terms of the outer variables, giving
\begin{equation}
g = \xi^{*2}\delta^2\left[\log\left(1/\epsilon\right) + \log(f)\right] + O(\epsilon).
\label{eq:gInnerOuter}
\end{equation}
To have an $O(1)$ term on the right of \eqref{eq:gInnerOuter} to match that on the left, we require $\log(1/\epsilon) = O(\delta^{-2})$.  Thus we let
\[
\log(1/\epsilon) = \ell_0\delta^{-2} + \ell_1\delta^{-1} + \ldots, \qquad \xi^* = \xi^*_0 + \delta \xi^*_1 + \ldots.
\]
Substituting these expansions into \eqref{eq:gInnerOuter}, along with those for $f$ and $g$ \eqref{eq:outerSolutionFirst}, \eqref{eq:outerSolutionNext}, we have from each order:
\begin{subequations}
\begin{align}
g_0 &= \xi_0^{*2} \ell_0 \\
g_1 &= \xi_0^{*2}\ell_1 + 2\xi_0^*\xi_1^* \ell_0 \\
g_2 &=  \xi_0^{*2}\log(f_0) + \mathrm{constant} \\
g_3 &= \xi_0^{*2}\frac{f_1}{f_0} + 2\xi_0^*\xi_1^* \log(f_0) + \mathrm{constant}
\end{align}
\label{eq:innerOuter}
\end{subequations}
Here the constant terms in the latter two equations involve higher order corrections ($\xi^*_2$, $\ell_2$, etc.).

The terms in the expansions for $\xi_0$ and $\log(1/\epsilon)$ are now found as follows.  The constants $\xi_0^*$ and $\xi_1^*$ are found by comparing the logarithmic terms in \eqref{eq:innerOuter}(c-d) to those in \eqref{eq:outerSolutionNext}.  The more singular term in $g_3$ $\propto f_0^{-1}$ then matches automatically.  The constants $\ell_0$ and $\ell_1$ are then found by comparing the constant terms in \eqref{eq:innerOuter}(a-b) to \eqref{eq:outerSolutionFirst}.  We thus find
\begin{equation}
\xi_0^* = a_{11}, \qquad \xi_1^* = a_{12}, \qquad \ell_0 = \frac{1}{a_{11}^2}, \qquad \ell_1 = -\left(\frac{1}{a_{11}^2} + \frac{2a_{12}}{a_{11}^3} \right),
\end{equation}
where $a_{11}$ and $a_{12}$ are the terms in the expansion of the slope at $\eta=0$ \eqref{eq:a1Terms}.

Finally, since $\epsilon$ is the independent parameter of the original problem, we invert the series expansions to find the front location as (logarithmic) series expansions in $\epsilon$.  Recall that in the original similarity coordinates, the apparent front is at $\eta^* = \delta \xi^*$.  Thus
\begin{subequations}
\begin{multline}
\delta \sim \ell_0^{1/2}\log(1/\epsilon)^{-1/2} + \frac{\ell_1}{2}\log(1/\epsilon)^{-1} \\
= \sqrt{\frac{\pi}{2}} \log(1/\epsilon)^{-1/2} - \frac{4\pi - 1 - \log 2}{6}\log(1/\epsilon)^{-1}
\label{eq:deltaFinalSeries}
\end{multline}
\begin{multline}
\eta^* \sim \xi_0^*\ell_0^{1/2} \log(1/\epsilon)^{-1/2} + \left(\frac{\xi_0^*\ell_1}{2} + \xi_1^*\ell_0\right) \log(1/\epsilon)^{-1} \\
= \log(1/\epsilon)^{-1/2} - \frac{\sqrt\pi}{2\sqrt 2} \log(1/\epsilon)^{-1}
\label{eq:etaSFinalSeries}
\end{multline}
\begin{multline}
a_1 \sim a_{11}\ell_0^{1/2} \log(1/\epsilon)^{-1/2} + \left(\frac{a_{11}\ell_1}{2} + a_{12}\ell_0\right) \log(1/\epsilon)^{-1} \\
= \log(1/\epsilon)^{-1/2} - \frac{\sqrt\pi}{2\sqrt 2} \log(1/\epsilon)^{-1}.
\label{eq:a1FinalSeries}
\end{multline}
\label{eq:finalSeries}%
\end{subequations}
Note that the expansions for the front location, $\eta^*$, and the slope at zero, $a_1$, are the same to this order.

\subsection{Comparison with numerical results}

We of course want to compare the above asymptotic result with the numerical results of Section \ref{sec:numerics}.  We plot the expansions \eqref{eq:finalSeries} in Figure~\ref{fig:epsilonCurvesEtaS} along with the numerical calculations of the same quantities.
The logarithmic dependence of the solution parameters on $\epsilon$ is certainly consistent with their observed numerical behaviour, even if the agreement is not very close at the level we can achieve numerically.
A major challenge in comparing the asymptotic solution directly to numerical results is the slow convergence of the logarithmic terms.  For $\epsilon = 10^{-6}$, which is around the smallest value we can attain with our numerical method, $[\log(1/\epsilon)]^{-1/2} \approx 0.27$, which is not very small.  At this value, the correction terms in the asymptotic formulas \eqref{eq:finalSeries} do improve the leading order result.  For $\epsilon = 10^{-6}$, the leading order asymptotic estimates of $\delta$ in \eqref{eq:deltaFinalSeries} and $\eta^*$ in \eqref{eq:etaSFinalSeries} are, respectively,
\[
\delta \sim \sqrt{\frac{\pi}{2}}[\log(1/\epsilon)]^{-1/2} \approx 0.34, \qquad \eta^* \sim [\log(1/\epsilon)]^{-1/2} \approx 0.27.
\]
Including the $O([\log(1/\epsilon)]^{-1})$ correction terms leads to estimates $\delta \approx 0.206$ and $\eta^* \approx 0.224$.  By comparison, the values estimated from the numerical solution in Section \ref{sec:numericalResultsForSmallEpsilon} are
\begin{equation}
\delta = 1-a_0 \approx 0.227, \qquad \eta^* \approx 0.241.
\label{eq:estimateParameters}
\end{equation}

\begin{figure}
\centering
\includegraphics[width=\textwidth]{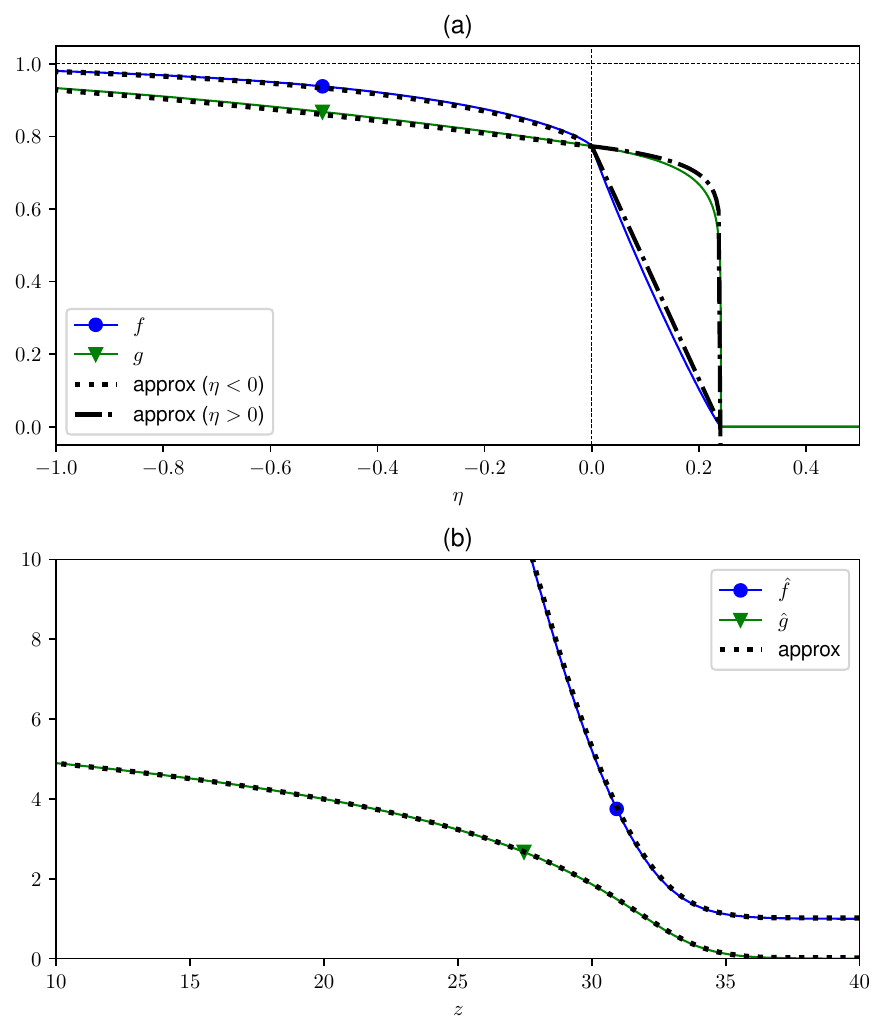}
\caption{(a) Comparison between the numerically computed similarity solution with $\epsilon=10^{-6}$, and approximate profiles given by \eqref{eq:negativeSolutionFirst} for $\eta < 0$ and \eqref{eq:outerNumericalFit}.  Due to the very slow convergence of the logarithmic terms as $\epsilon \to 0$, numerical estimates of $\delta$ and $\eta^*$ are used in fitting these profiles.  (b) The asymptotic and numerical profiles in the inner region near $\eta^*$.  Since the asymptotic solution is accurate to $O(\epsilon)$ in this region, the agreement is excellent.}
\label{fig:epsilonNearZero}
\end{figure}

However, by using the numerical estimates for $\delta$ and $\eta^*$ in \eqref{eq:estimateParameters} instead of asymptotic approximations, we can find encouraging agreement with the solution profiles in each region.  In Figure~\ref{fig:epsilonNearZero}(a) we plot the asymptotic solution for $\eta < 0$ and $0 < \eta < \eta_0$.  Here we use the $O(\delta)$ expansions \eqref{eq:negativeSolutionFirst} for $\eta < 0$, and for $\eta > 0$ use the approximations for $f$ and $g$ that satisfy the conditions at $\eta=0$ and include their first nonuniform terms, that is:
\begin{equation}
g \approx 1 - \delta + \delta^2\sqrt{\frac{2}{\pi}}\log\left(1-\frac{\eta}{\eta^*}\right), \qquad f \approx (1-\delta)\left(1- \frac{\eta}{\eta^*}\right).
\label{eq:outerNumericalFit}
\end{equation}

In Figure~\ref{fig:epsilonNearZero}(b) we compare the solution close to $\eta^*$, scaled into variables $z$, $\hat f$ and $\hat g$, with the inner solution \eqref{eq:innerSolution}.  To do so we also have to numerically estimate the translational constant $z_0$ (this constant is exponentially small in $\delta$, so in terms of the asymptotic analysis cannot not be estimated from matching to any finite number of terms of the outer problem).  Otherwise, the agreement is practically perfect, which is to be expected given the error in the inner region is of order $\epsilon$, rather than $\delta$.

\section{Discussion}
\label{sec:discussion}

Our main finding in this study is the singular nature of the limit of a system that couples degenerate diffusion and degradation of surrounding environment, as we attempt to approach a system that could exhibit a moving interface (that is, a point where the diffusivity vanishes).  The fact that the speed of propagation vanishes in this limit, and gradients of the solution become infinite, is strong evidence that the system cannot exhibit such a moving interface.

The analysis undertaken for the idealised system \eqref{eq:pde} is also relevant for more complex models of the form \eqref{eq:pdeWithReactions} for which coupled degenerate diffusion and degradation is the driving factor.
Consider a model of the form \eqref{eq:pdeWithReactions} with $F(\tilde u)$ such that $F(1) = 0$, starting at a piecewise-constant initial condition such as \eqref{eq:pdeIC}.  For early times, we may use the same similarity scaling as \eqref{eq:similarityAnsatz}, which will then be the leading order behaviour of solutions for small $t$, that is, 
\[
\tilde u = 1-tf(\eta) + O(t), \qquad v(x,t) = g(\eta) + O(t), \qquad \eta = \frac{x-x^*_0}{t}.
\]
The effect of the reaction terms, as well as the diminishing decay of $\tilde u$ as $\tilde u$ itself decreases, induces an $O(t)$ correction to the leading order solution, which is the same similarity solution as considered in this study.  More generally, we expect the inner problem \eqref{eq:innerProblem2} to be the relevant inner problem near any apparent moving interface (not just in the context of a similarity solution), so that the logarithmic dependence of speed on a regularisation of order $\epsilon$ near the apparent front will hold in other contexts, for example travelling wave solutions.  We also note that for early times, the same similarity solution will hold near a step discontinuity in versions of \eqref{eq:pdeWithReactions} posed in higher spatial dimensions.

The logarithmic dependence of speed on the small parameter $\epsilon$ explains previous observations made chiefly from PDE simulations in the literature~\cite{colson2021travelling, el2021travelling, mascia2024numerical}.
It is not surprising that numerical methods applied to systems of the form \eqref{eq:pdeWithReactions} with an initially compactly supported condition will appear to evolve in time with a well-defined front.  In a numerical scheme, $\epsilon$ is representative of the mesh size of the spatial discretisation.  It is very challenging in a discretised PDE simulation to obtain a mesh size less than $\sim 10^{-3}$ or so, and for such a value of $\epsilon$, the apparent speed of propagation of a front (including, for example, the apparent speed of a travelling wave) will appear to be a positive quantity.  The careful examination of the convergence of numerical simulations, as was done in \cite{mascia2024numerical} for example, is thus crucial for these forms of degenerate diffusion systems.

While we have considered the regularising parameter $\epsilon$ as the height in front of an evolving profile, it could be considered equivalent to a number of regularisations that will produce similar results.  
For example, in a continuum limit of a discrete process (e.g. an agent based model of individual cells), the (nondimensional) individual cell size could be considered a small but nonzero regularising parameter.  In each case, the details of the $O(\epsilon)$ inner problem near an apparent interface would be different, but the matching with the outer problem would result in the same logarithmic dependence.

The regularising effect of height $\epsilon$ near a moving front has interesting parallels to the contact-line paradox in interfacial fluid dynamics, in which the classical laws of hydrodynamics (in particular the no-slip boundary condition) cannot model the motion of contact line in the spreading of a droplet, for instance \cite{huh1971hydrodynamic}.  In such cases, a number of regularisations (such as disjoining pressure or Navier slip) are posited that regularise the singularity at the contact line, and the contact line speed is logarithmically dependent on the small parameter in the regularisation \cite{bonn2009wetting,snoeijer2013moving}.

Another example of logarithmic dependence on a small regularisation from mathematical biology is the effect of perturbing the reaction term in the Fisher--KPP equation on its travelling wave speed, for example the cut-off effect considered in \cite{brunet1997shift,dumortier2007critical}.  Here the perturbation is introduced to take into account the discrete size of cells.  Although the limit is not singular in this case (the standard minimum Fisher--KPP travelling wave is recovered in the limit), the effect on the wave speed is logarithmic, and so can deviate significantly from that expected.

We now consider some possible extensions of this work.  The similar class of cross-diffusion models considered \cite{crossley2025existence,crossley2023traveling,simpson2024discrete}, which may, at least in an idealised form, be represented by the system
\[
\pfrac{u}{t} = v, \qquad \pfrac{v}{t} = \pfrac{}{x}\left(u\pfrac{v}{x} + \phi(v)\pfrac{u}{x}\right),
\]
for some function $\phi(v)$.  For an initially piecewise condition such as \eqref{eq:pdeIC}, a system of this form will also be amenable to a similarity scaling similar to the one carried out in this study. We foresee such systems having the same logarithmic dependence on any regularisation near an apparent interface, although the details of the similarity solution will depend on the form of $\phi$.

Finally, if one wanted a system similar to \eqref{eq:pde} that allowed for actual moving interfaces, one approach could be to weaken the degeneracy of the diffusivity by introducing a power-law relationship:
\begin{equation}
\pfrac{u}{t} = v, \qquad \pfrac{v}{t} = \pfrac{}{x}\left(u^m\pfrac{v}{x}\right),
\label{eq:pdeWithPower}
\end{equation}
where the exponent $m \in (0,1)$.  This alteration is likely to allow for (weak) solutions with moving interfaces, as if we expand solutions of \eqref{eq:pdeWithPower} close to an interface $x^*(t)$ where $u$ and $v$ both vanish, we find a valid balance:
\[
u \sim U_0(x^*-x)^{1/m}, \qquad v \sim \frac{1-m}{m^2}U_0^{m+1}(x^*-x)^{1/m-1}, \qquad x \to x^{*-}
\]
where the velocity of the interface $\dot x_0^* = (1-m)A^m/m$.  Of course, these terms are clearly singular in the limit $m\to 1^-$ as we approach the original system \eqref{eq:pde}.  We will consider such problems in future studies.

\end{document}